\documentclass[10pt]{amsart} \usepackage{amssymb, verbatim, epic, eepic, epsfig}

\def\dist{{\hbox{\rm dist}}}

\def\R{{{\mathbf{R}}}}

\def\T{{\hbox{\bf T}}}

\def\emph#1{{\it #1}} \def\textbf#1{{\bf #1}}

\def\pmax{p_{\max}}
\def\pmin{p_{\min}}

\def\vol{\operatorname{vol}}
\def\cH{{\mathcal H}} 
\def\Ricci{\operatorname{Ricci}}

\def\RR{\mathbb{R}}

\def\ang{{\operatorname{ang}}}

\def\sc{\operatorname{sc}}
 
\def\Id{\operatorname{Id}}
\def\ang#1{\langle #1 \rangle}

\def\sc{\operatorname{sc}}
\def\MMsc{\overline{M}^2_{\sc}}
\def\Diagsc{\Delta_{\sc}}
\def\Diagscsp{\Delta_{\sc, k}}

\def\rb{\operatorname{rb}}
\def\rbz{\operatorname{rb0}}
\def\lb{\operatorname{lb}}
\def\lbz{\operatorname{lb0}}
\def\bfc{\operatorname{bf}}
\def\bfacez{\operatorname{bf0}}
\def\sf{\operatorname{sf}}
\def\ff{\operatorname{ff}}
\def\bbdy{\operatorname{bb}}
\def\zf{\operatorname{zf}}
\def\MMscsp{\overline{M}^2_{\operatorname{sc}, k}}

\def\Msp{\overline{M}_{k}}
\def\interior{\operatorname{int}}

\def\Real{\operatorname{Re}}
\def\Imag{\operatorname{Im}}

\newcommand{\Mbar}{\overline{M}}

\newcommand{\CI}{{\mathcal{C}^\infty}}

\renewcommand\lb{\operatorname{lb}}
\renewcommand\rb{\operatorname{rb}}

\newcommand\CIdot{\dot C^\infty}

\theoremstyle{plain} \newtheorem{theorem}[subsection]{Theorem}
  
  \newtheorem{proposition}[subsection]{Proposition}
  \newtheorem{lemma}[subsection]{Lemma}
  \newtheorem{corollary}[subsection]{Corollary}
\theoremstyle{remark} 
  \newtheorem{rem}[subsection]{Remark}
  
  \newtheorem{openproblem}[subsection]{Open Problem}

\theoremstyle{definition} 
  
\numberwithin{equation}{section}

\include{psfig}

\begin{document}

\title[Riesz transform for manifolds with Euclidean ends]{Riesz transform and $L^p$ cohomology for manifolds with Euclidean ends}
\author{Gilles Carron}\address{D\'epartement de Math\'ematiques, Universit\'e de Nantes,
44322,  Nantes, FRANCE} \email{gilles.carron@math.univ-nantes.fr}
\author{Thierry Coulhon}\thanks{T.C. is supported in part  by the European Commission (IHP
Network ``Harmonic
Analysis and Related Problems'' 2002-2006, Contract HPRN-CT-2001-00273-HARP). He also acknowledges for this work the support of Macquarie University and Australian National University}
\address{D\'epartement de Math\'ematiques, Universit\'e de Cergy-Pontoise,
95302,  Pontoise, FRANCE} \email{thierry.coulhon@math.u-cergy.fr}
\author{Andrew Hassell} \thanks{A.H.\ is supported in part by an Australian
Research Council Fellowship and acknowledges the support of Universit\'e de Nantes.}  \address{Department of Mathematics, ANU,
Canberra, ACT 0200, AUSTRALIA} \email{hassell@maths.anu.edu.au}



\vspace{-0.3in}
\begin{abstract}Let $M$ be a smooth Riemannian manifold which is the union of a compact part and a finite number of Euclidean ends, $\RR^n \setminus B(0,R)$ for some $R > 0$, each of which carries the standard metric. Our main result is that the Riesz transform on $M$ is bounded from $L^p(M) \to L^p(M; T^*M)$ for $1 < p < n$ and unbounded for $p \geq n$ if there is more than one end. It follows from known
results that in such a case the Riesz transform
on $M$ is bounded for $1 < p \leq 2$ and unbounded for $p > n$; the result is new for $2 < p \leq n$. We also give some heat kernel estimates on such manifolds. 

We then consider the implications of boundedness of the Riesz transform in $L^p$ for some $p > 2$ for a more general class of manifolds. Assume that $M$ is a $n$-dimensional complete manifold satisfying the Nash inequality and with an $O(r^n)$ upper bound on the volume growth of geodesic balls. We show that boundedness of the Riesz transform on $L^p$ for some $p > 2$ implies a Hodge-de Rham interpretation of the $L^p$ cohomology in degree $1$, and that the map from $L^2$ to $L^p$ cohomology in this degree  is injective.
\end{abstract}

\maketitle


\section{Introduction}\label{sec:introduction}
Let $M$ be a complete Riemannian manifold with infinite measure.
The Riesz transform $T$ on $M$ is the operator 
\begin{equation}
f \to d \Delta^{-1/2} f ,
\end{equation}
where $\Delta$ is the positive Laplace operator on $M$. The Riesz transform is always a bounded map from $L^2(M)$ to $L^2(M; T^*M)$. It is of interest to figure out the range of $p$ for which $T$ extends to a bounded map $L^p(M) \to L^p(M; T^*M)$. Equivalently, we can ask whether
\begin{equation*}
\| \, |d f | \, \|_p \leq \| \Delta^{1/2} f \|_p \text{ for all } f \in C_c^\infty(M).
\end{equation*}It has been shown in \cite{CD} that the Riesz transform is bounded on $L^p$ for $1 < p < 2$ as soon as the manifold satisfies 
the doubling property as well as a natural heat kernel upper bound. The situation for $p > 2$ is more complicated:  there is some understanding of what happens in the more restricted class of manifolds
satisfying upper and lower Gaussian estimates for the heat kernel (see \cite{ACDH}, \cite{AC}), and  it is also known that the above more general assumptions
do not  imply the boundedness of the Riesz transform for all $p>2$. The counterexample is simply the connected sum of two copies 
of $\RR^n$, where one easily sees that the Riesz transform is unbounded for  $p>n$ (see \cite{CD}), precisely because the  boundedness for such $p$ would imply a lower Gaussian heat kernel estimate on this manifold, which is known to be false. The first aim of the present article is to find out what happens in the remaining range
$2<p\le n$, and to treat more generally by the same token manifolds with a finite number of Euclidean ends.
Our second aim is to give, in a somewhat more general class of manifolds,  a cohomological consequence of the boundedness of the Riesz transform for some $p>2$, which explains the unboundedness in the range $p\ge n$ in the examples just mentioned. 

\bigskip

On $\RR^n$, it is a classical result in harmonic analysis that the Riesz transform is bounded on $L^p$ for all $1 < p < \infty$. In this case, the kernel of the Riesz transform is given by
\begin{equation*}
c d \frac1{|z-z'|^{n-1}} = c' \frac{(z-z')_i dz^i}{|z-z'|^{n+1}}.
\end{equation*}
(We shall use $z$ for a Euclidean coordinate throughout this paper, while the prime denotes the `right variable' of a kernel $K(z,z')$; also, $\hat z$ denotes $\frac{z}{|z|}$ and we write $r, r'$ for $|z|, |z'|$.)
It is useful to compactify $\RR^n$ to a compact manifold with boundary $Z$ by adding a sphere at infinity, and using $|z|^{-1}$ as a boundary defining function\footnote{This means that smooth functions near the boundary of the compactification are given precisely by smooth functions of $\hat z$ and $1/|z|$.} for it. If we consider this kernel of $\Delta^{-1/2}$ near the `right boundary' of $Z \times Z$, namely $Z \times \partial Z$, where $|z'| \to \infty$, we see that it has an expansion given by
\begin{equation*}
\frac1{|z-z'|^{n-1}} = \frac1{|z'|^{n-1}} \big( 1 - 2 \frac{z \cdot \hat z'}{|z'|} + \frac{|z|^2}{|z'|^2} \big)^{-\frac{n-1}{2}} \sim \frac1{|z'|^{n-1}} + (n-1) \frac{ z \cdot \hat z'}{|z'|^n} + \dots
\end{equation*}
The leading power in this expansion is $|z'|^{-n+1}$ and the coefficient multiplying it is $1$. This kernel, by itself, is not bounded on $L^p$ for $p \geq n$ because the decay of the kernel puts it only in $L^{n/(n-1) + \epsilon}$, $\epsilon > 0$, as a function of $z'$ with $z$ fixed,  so it can only be boundedly paired with elements in $L^p$ for $p < n$. However, when we apply a $z$-derivative, the leading term is killed (since $d 1 = 0$) and the kernel of $T$ decays at one order better, namely $|z'|^{-n}$. This allows pairing with elements of $L^p$ for any $p < \infty$.
Let $M$ be a manifold with Euclidean ends, and assume the number of ends is at least two.

It is now relatively easy to explain why the Riesz transform on $M$ is not bounded for $p \geq n$. For simplicity we shall assume here that $M$ has exactly two Euclidean ends. Let us compactify $M$ to a compact manifold $\Mbar$ in the analogous way to $Z$ above. The boundary of $\Mbar$ is then the disjoint union of two $n-1$-spheres, which we shall denote $\partial \Mbar_+$ and $\partial \Mbar_-$.
It turns out that the kernel of $\Delta^{-1/2}$ has a similar expansion at $\Mbar \times \partial \Mbar$, of the form
\begin{equation} 
\sum_{j=n-1}^\infty |z'|^{-j} a_j(z), \quad a_j \in \CI(\Mbar \times \partial M),
\end{equation}
but it is no longer true that the leading term $a_{n-1} = 1$. Rather, $a_{n-1}$ is the harmonic function on $\Mbar$ that equals $1$ on $\partial \Mbar_+$, and zero on $\partial \Mbar_-$. As a consequence, applying the derivative operator in the left variable $z$ does not make this leading term disappear, since $d a_{n-1} \neq 0$. Hence the kernel of $T$ only decays to order $n-1$ at the right boundary of $\Mbar \times \Mbar$, and therefore can only be paired boundedly with elements of $L^p$ for $p < n$.

In this paper we shall prove

\begin{theorem}\label{main} Let $M$ be a complete $\CI$ Riemannian manifold of dimension $n \geq 3$ which is the union of a compact part and a finite number of Euclidean ends. Then the Riesz transform is bounded from $L^p(M)$ to $L^p(M; T^*M)$ for $1 < p < n$, and is unbounded on $L^p$ for all other values of $p$
if the number of ends is at least two. 
\end{theorem}
Our method is to analyze the kernel $\Delta^{-1/2}$ based on the formula
\begin{equation}
\Delta^{-1/2} = \frac2{\pi} \int_0^\infty (\Delta + k^2)^{-1} \, dk.
\label{sp-int}\end{equation}
Since $M$ is a manifold which is conic at infinity, the Laplacian $\Delta$ on $M$ lies in the class of scattering differential operators \cite{scatmet} and we can use methods from the scattering calculus to analyze the kernel of $\Delta^{-1/2}$. We shall analyze the kernel of $\Delta^{-1/2}$ rather precisely and work out the leading term in the expansion at the right boundary. F{}rom this, it will be straightforward to analyze the kernel of $T$ and to prove the theorem. The plan of the paper is as follows. We briefly describe the scattering calculus in section~\ref{scatt}. In sections~\ref{para} -- \ref{RT} we prove the theorem using the analysis of the resolvent of the Laplacian on asymptotically conic spaces in \cite{HV} as a model. We give some large time asymptotics on derivatives of the heat kernel on manifolds with Euclidean ends in section~\ref{heat}. 

In section~\ref{coh} we change point of view and consider a much more general class of manifolds, namely complete manifolds $M$ of dimension $n$ satisfying the Nash inequality and with a uniform upper bound $O(r^n)$ on the volume of geodesic balls of radius $r$. We \emph{assume} that the Riesz transform on $M$ is bounded on $L^p$ for some $p > 2$ and give several geometric and topological consequences: a Hodge-de Rham interpretation of the $L^p$ cohomology of $M$ (Proposition~\ref{HdR}), injectivity of the map from $L^2$ to $L^p$ cohomology (Lemma~\ref{inj}) and derive a contradiction if $p > n$ and $M$ has at least two ends (Corollary~\ref{nb}), thus generalizing the unboundedness part of Theorem~\ref{main} for $p > n$ to this larger class of manifolds. In the final section we discuss our results in the context of previously known examples and pose some open problems.


\section{Scattering Calculus}\label{scatt}
As noted above, we shall use the scattering calculus \cite{scatmet} to analyze the kernel of $\Delta^{-1/2}$ on manifolds with several Euclidean ends. The scattering calculus is expressed in terms of compactifications on $M$ and, especially, of the double space $M^2$ which carries the kernel of the resolvent and of the operator $\Delta^{-1/2}$. The space $M$ is compactified by adding a sphere at infinity $S^{n-1}$ for each Euclidean end, and declaring $r^{-1}= 1/|z|$ and $\hat z = z/|z|$ to be  local coordinates near a boundary point; in particular, $r^{-1}$ is taken to be a defining function for the boundary (which we shall sometimes refer to as `infinity'). We sometimes use $x = r^{-1}$ to denote this boundary defining function and $y = \hat z$ to denote boundary coordinates, extending to a collar neighbourhood of the boundary, as is customary when using the scattering calculus. The metric takes the form
\begin{equation*}
dr^2 + r^2 h(y, dy) = \frac{dx^2}{x^4} + \frac{h(y, dy)}{x^2}
\end{equation*}
at each end, and is therefore a \emph{scattering metric} as defined in \cite{scatmet} (of a particularly simple form, being an exact conic metric near infinity).

We are mostly interested in the case when the number of ends is at least two. In fact, for the sake of clear exposition we shall assume from now on that the number of ends is exactly two, although all proofs in this paper generalize in an obvious way to any finite number of ends. We shall label these ends $+$ and $-$, thus for example we shall use $z_+$ as the Euclidean variable on the positive end, and $z_-$ for the variable on the negative end; when it is not necessary to stipulate which end is being considered, we shall just use $z$.

It is not so obvious which compactification of $M^2$ is most appropriate for dealing with the Schwartz kernels of operators such as the Laplacian, or functions of the Laplacian, on manifolds with Euclidean ends. There are several different asymptotic regimes of interest when dealing with such kernels.  One regime, the `near-diagonal' regime, is when the two variables $z, z'$ remain a finite distance apart as they both go to infinity. Another is when they both go to infinity with the ratio $r/r'$ approaching a limit and with $\hat z,\hat z'$ both approaching a limit. Finally there is the case that one variable approaches a limit, while the other remains fixed. The kernel has different behaviour in each of these asymptotic regimes, so they need to be represented by distinct parts of the boundary of the compactification.
The space $(\Mbar)^2$ is thus  too `small' a compactification of $M^2$ for our purposes, since it only has the third regime distinguished; the first two are squashed into the corner.

It turns out that there is a space denoted $\MMsc$, the `scattering double space', which satisfies these criteria. It is obtained by performing two blowups on $M^2$. The first is blowing up the corner $(\partial \Mbar)^2$, creating the so-called b-double space, and the second is blowing up the boundary of the diagonal (which lifts to the b-double space to be transverse to the boundary, hence this blowup is well-defined). Each asymptotic regime is represented by a boundary hypersurface of $\MMsc$. The first, `near-diagonal' regime is represented by the boundary hypersuface created by the second blowup, denoted $\sf$ for `scattering face'; the second is represented by the boundary hypersurface created by the first blowup, denoted $\bfc$ for `b-face' (since it is present in the b-calculus) and the third regime is represented by the two boundary hypersurfaces $\partial \Mbar \times \Mbar$, $\Mbar \times \partial \Mbar$ of $(\Mbar)^2$, denoted $\lb$ and $\rb$ for `left boundary' a
 nd `right boundary'. Note that when $M$ has $k$ ends, then $\bfc$ has $k^2$ components and $\sf$, $\lb$ and $\rb$ each have $k$ components. \begin{figure}\centering
\epsfig{file=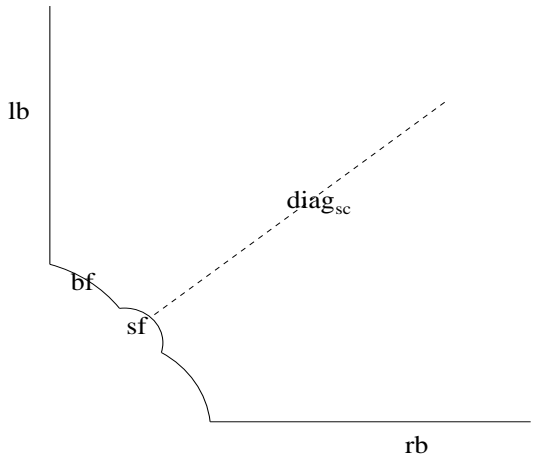,width=7cm,height=6cm}
\caption{The double scattering space}
\label{fig:MMsc}
\end{figure}

The structure of $\sf$ and $\bfc$ is as follows. Each component of $\sf$ is naturally diffeomorphic to $\overline{\RR^n} \times S^{n-1}$, and $z - z'$ and $\hat z$
are coordinates on the interior of each component of $\sf$. Each component of $\bfc$ is naturally isomorphic to a blowup of the space $S^{n-1} \times S^{n-1} \times [0,1]$, with coordinates $(\hat z, \hat z', |z|(|z| + |z'|)^{-1})$; the blowup is of the submanifold $\{ \hat z = \hat z', |z|(|z| + |z'|)^{-1} = 1/2 \}$ which corresponds to the boundary of the diagonal.

The \emph{scattering calculus} is an algebra of pseudodifferential operators on $M$ which is defined by the properties of their Schwartz kernels. Namely, $A$ is a scattering psuedodifferential operator of order $(m,0)$ on $M$ iff the kernel of $A$, when lifted to $\MMsc$, is conormal\footnote{In other words, the kernel of $A$ has a singularity at the diagonal characteristic of pseudodifferential operators of order $m$, and this holds smoothly up to the boundary,  in the sense that it could be extended across the boundary as a conormal distribution. See \cite{Ho}, section 18.2, for the precise definition.} of order $m$ at the diagonal of $\MMsc$ smoothly up to the boundary $\sf$, is smooth elsewhere at $\sf$, and is rapidly decreasing at $\bfc$, $\lb$ and $\rb$. The resolvent of the Laplacian $(\Delta - \lambda^2)^{-1}$ is a scattering pseudodifferential operator of order $(-2, 0)$ on $M$ for $\Real \lambda \neq 0$. In fact, the structure of the resolvent on the spectrum, i.e.
  the kernel of $(\Delta - (\lambda \pm i0)^{2})^{-1}$ for real $\lambda$, can also be described on $\MMsc$, although here the kernel is no longer rapidly decreasing at $\bfc$, $\lb$ and $\rb$, rather it is a `Legendrian distribution' \cite{HV}. Our approach is partly modelled on the analysis in this paper. However, we take advantage of the assumption here that $M$ has exact Euclidean ends, which leads to great simplifications over the analysis of \cite{HV} since we can exploit the well-known explicit formulae for the resolvent of the Laplacian on $\RR^n$ and use these as ingredients for a parametrix of the resolvent kernel on $M$, thereby avoiding the need to use Legendrian distributions in this paper.

\emph{Notation.} We write $z$ for a Euclidean variable $z = (z_1, z_2, \dots, z_n) \in \RR^n)$ and write $\ang{z} = \sqrt{1 + |z|^2}$, while $x$ is used for  $|z|^{-1}$, or, sometimes, where more convenient, for $\ang{z}^{-1}$. For a manifold with corners $X$, we write $\CIdot(X)$ for the space of smooth functions which vanish to infinite order at the boundary of $X$. We use notation $[X; S_1, S_2, \dots S_n]$ to denote the blowup of $X$ at the submanifolds $S_1$, $S_2$, \dots (in that order).

\section{Parametrix construction}\label{para}

To analyze the operator $\Delta^{-1/2}$ we return to the formula \eqref{sp-int}.
We first observe that the off-diagonal terms in the kernel of $\Delta^{-1/2}$ come from a neighbourhood of zero in the integral \eqref{sp-int}. Indeed, let $s_0(k)$ be a cutoff function equal to $1$ in a neighbhourhood of $k=0$ and equal to zero outside a compact set, and let $s_1(k) = 1 - s_0(k)$. Then we may insert the factor $1 = s_0(k) +s_1(k)$ into the integral \eqref{sp-int}. With the factor $s_i$ inserted, the integral gives a function $g_i(\Delta)$ of $\Delta$, $i = 0$ or $1$, where
$$
g_1(t) = \int_0^\infty s_1(k) \frac1{k^2 + t^2} \, dk 
$$
is easily checked to be a classical symbol of order $-1$. By the symbolic functional calculus \cite{HV:Symbolic}, this term is a scattering pseudodifferential operator of order $-1$, hence $d g_1(\Delta)$ is a scattering pseudodifferential operator of order zero. It is therefore bounded on $L^p$ for all $1 < p < \infty$ \cite{Stein}. So we are reduced to studying $d g_0(\Delta)$, given by  the integral \eqref{sp-int} with factor $s_0(k)$ inserted.

We shall write down a fairly explicit parametrix for $(\Delta + k^2)^{-1}$ for small $k$. In doing so, we need to consider the different asymptotics that this kernel takes when $k= 0$ and $k \neq 0$. Indeed, on $\RR^n$ the kernel decays as $|z-z'|^{-(n-1)/2}$ for $k \neq 0$ and $|z-z'|^{-n+2}$ for $k = 0$, as $|z-z'| \to \infty$, which (except when $n=3$) is a different rate. This can be encoded geometrically by blowing up at the boundary when $k = 0$. Consider the space
\begin{equation}
\MMscsp = [\MMsc \times [0, k_0]; \bfc \times \{ 0 \}; \lb \times \{ 0 \}; \rb \times \{ 0 \}].
\label{scsp-blowup}\end{equation}
We shall denote the boundary hypersurfaces which are the lifts of $\bfc \times [0, k_0]$, $\lb \times [0, k_0]$ and $\rb \times [0, k_0]$ to $\MMscsp$ by $\bfc$, $\lb$ and $\rb$, and $\sf \times [0, k_0]$ by $\sf$; this is of course an abuse of notation, but in the context it will always be clear whether it is a boundary hypersurface of $\MMsc$ or $\MMscsp$ that is referred to.
We shall denote the new boundary hypersurfaces corresponding to the three blowups by
$\bfacez$, $\lbz$ and $\rbz$, according as they arise from the first, second or third blowups in \eqref{scsp-blowup} respectively, and we shall denote $\MMsc \times \{ 0 \}$ by $\zf$, for `zero face'.
We also define $\Diagscsp$ to be $\Diagsc \times [0, k_0] \subset \MMscsp$. Let $\chi$ be a smooth function on $\MMscsp$ which is equal to one in a neighbourhood of $\Diagscsp$, and whose support meets the boundary of $\MMscsp$ only at $\sf$ and $\zf$.
\begin{figure}\centering
\epsfig{file=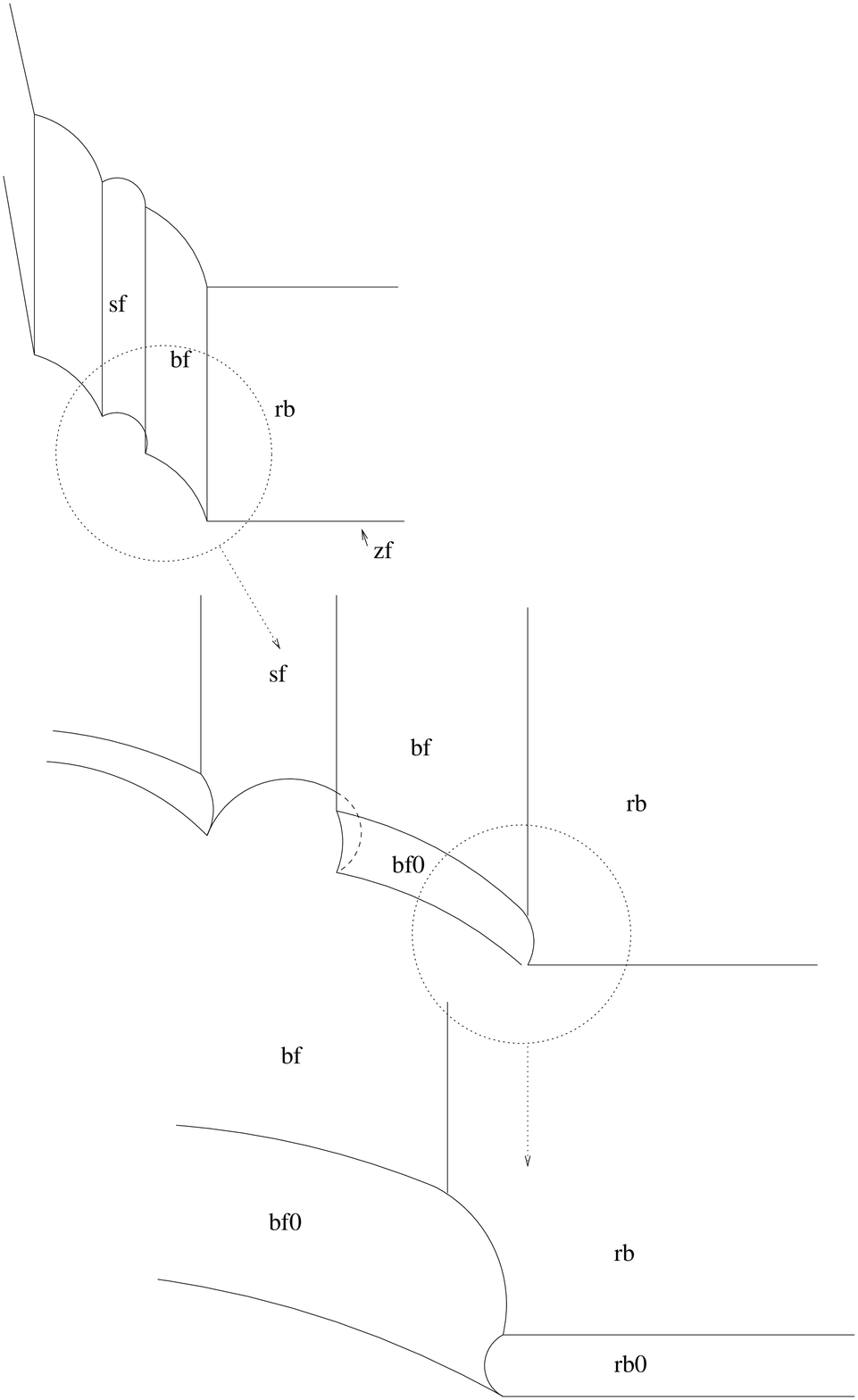,width=9cm,height=15cm}
\caption{Blowing up $\MMsc \times [0, k_0]$ to produce $\MMscsp$}
\label{fig:MMscsp}
\end{figure}

We recall the well-known expression for the resolvent kernel $(\Delta + k^2)^{-1}$ on $\RR^n$ 
for $n \geq 3$:
\begin{equation}(\Delta + k^2)^{-1} = \frac{e^{-k|z-z'|}}{|z-z'|^{-n+2}} f_n(k|z-z'|),
\label{res}\end{equation}
where $f_n(t)$ is symbolic of order $(n-3)/2$ as $t \to \infty$, while it is $O(1)$ and has a classical expansion in powers and logarithms as $t \to 0$. In fact, $f_n$ is a polynomial of order $(n-3)/2$ when $n \geq 3$ is odd.
It is straightforward to check

\begin{lemma}\label{conormal} Let $Z$ be the compactification of $\RR^n$. Then the resolvent kernel $(\Delta + k^2)^{-1}$ is such that $\chi (\Delta + k^2)^{-1}$ is conormal at $\Diagscsp$, and
\begin{equation}
(1 - \chi) (\Delta + k^2)^{-1} \in \rho_{\sf}^0  (\rho_{\bfc} \rho_{\lb} \rho_{\rb})^{\infty}(\rho_{\bfacez}\rho_{\lbz}\rho_{\rbz})^{n-2} \CI(\MMscsp).
\label{conorm}\end{equation}
Here, `conormal to $\Diagscsp$' means that the kernel is conormal  in $z-z'$ which defines $\Diagscsp$ and smooth in the remaining variables $x = |z|^{-1}, \hat z, k$, uniformly up to the boundary.
\end{lemma}
For example, let us check the statement of the lemma near the triple intersection 
$\bfacez \cap \rbz \cap \rb$. Coordinates near this codimension three corner are $\hat z, \hat z'$
 and boundary defining functions $\rho_{\rbz} = k/x = k|z|$ for $\rbz$, $\rho_{\rb} = x'/k = 1/(k|z'|)$
  for $\rb$ and $\rho_{\bfacez} = x$ for $\bfacez$. Near this corner, $|z'|$ is much larger than $|z|$
   so we may expand
   $$k|z-z'| = k|z'| \big( 1 - \frac{2z \cdot \hat z'}{|z'|} + \frac{|z|^2}{|z'|^2} \big)^{1/2} =
    \frac1{\rho_{\rb}} \big( 1 - 2\hat z \cdot \hat z' \rho_{\rb} \rho_{\rbz} + (\rho_{\rb} \rho_{\rbz} )^2 \big)^{1/2},$$from which it is easy to check that \eqref{conorm} holds.

We also need a single space version of this space. Let 
\begin{equation}
\Msp = [\Mbar \times [0, k_0]; \partial M \times \{ 0 \}].
\end{equation}
Denote the boundary hypersurfaces $\bbdy$, $\zf$ and $\ff$ which arise from $\partial \Mbar \times [0, k_0]$, $\Mbar \times \{ 0 \}$ and from the blowup, respectively, and denote corresponding boundary defining functions by $\rho_{\bbdy}$, $\rho_{\zf}$ and $\rho_{\ff}$. Again, it will always be clear in context whether $\zf$ refers to the zero-face of $\Msp$ or $\MMscsp$. \begin{figure}\centering
\epsfig{file=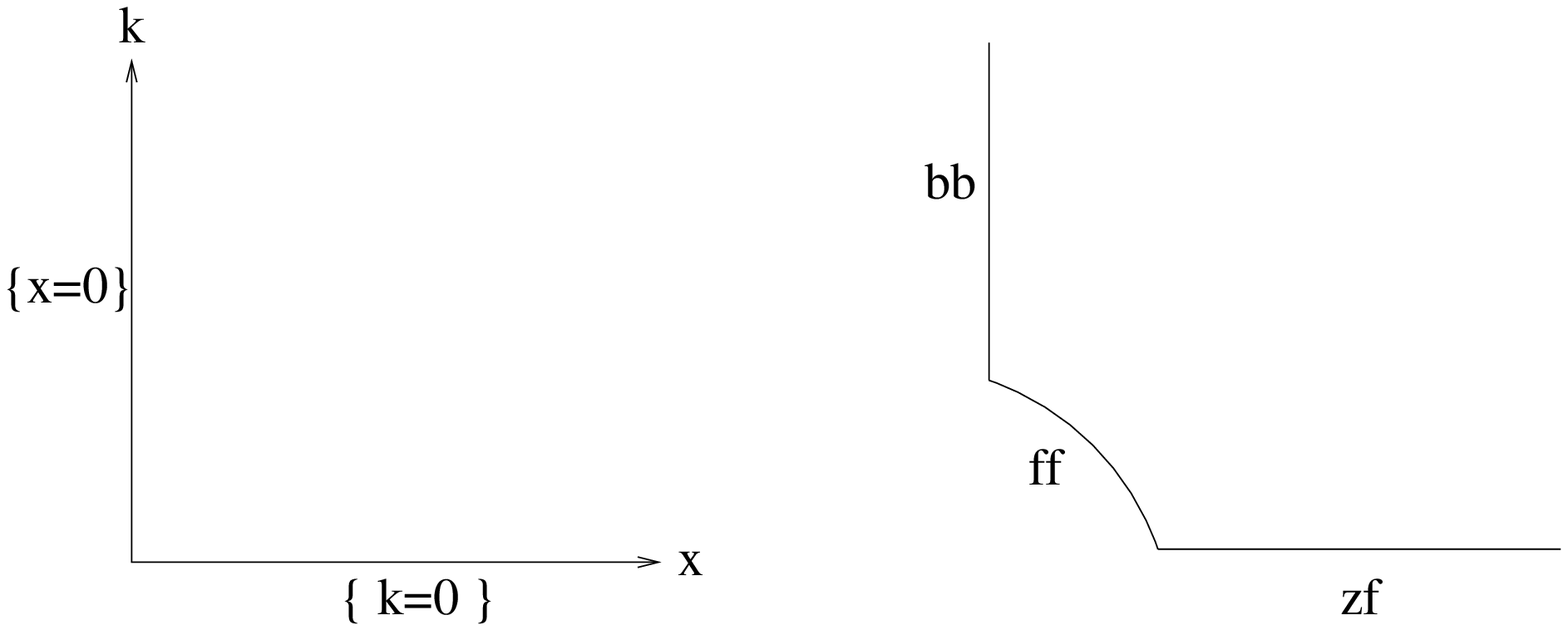,width=9cm,height=4cm}
\caption{Blowing up $M \times [0, k_0]$ to produce $\Msp$}
\label{fig:Msp}
\end{figure}We have 

\begin{lemma}\label{const} Let $v \in \CIdot(\Mbar)$. Then there is a function $u \in \rho_{\bbdy}^{\infty} \rho_{\ff}^{n-2} \CI(\Msp)$, such that $(\Delta + k^2) u | \zf$ is equal to $v$, and $(\Delta + k^2) u$ vanishes to infinite order at both $\ff$ and $\bbdy$.
\end{lemma}

\begin{proof}We first use results from \cite{tapsit} to show that we can solve $\Delta f = v$ on $M$. The Laplacian on an asymptotically Euclidean manifold $M$ may be written in the form $$\Delta = x^{n/2+1} P x^{n/2 - 1},$$where $P$ is an elliptic b-differential operator on $M$. A short computation shows that near infinity, $P$ takes the form $$-(x \partial_x)^2 + \big( \frac{n-2}{2} \big)^2 + x^2 \Delta_{S^{n-1}}, \quad x = |z|^{-1},$$where $\Delta_{S^{n-1}}$ is the standard Laplacian on the $(n-1)$-sphere. This is a strictly positive operator, so $P$ is `totally elliptic', and hence is Fredholm acting between the b-Sobolev spaces\footnote{Here $L^2_b(M)$ is the $L^2$-space with respect to the b-metric $g_b = x^2 g$; thus $L^2_b(M) = \ang{z}^{n/2} L^2(M)$. Also $H^2_b(M)$ is the b-Sobolev space of order $2$, defined as the set of functions $g \in L^2(M)$ such that $Qg \in L^2(M)$ for all b-differential operators of order $2$ on $M$ (near infinity, such  operators take the f
 orm $\ang{z}^2 \sum_{i,j} a_{ij} \partial_{z_i} \partial_{z_j}$, with $a_{ij} \in \CI(\Mbar)$).} $H^{2}_b(M) \to L^2_b(M)$.  Thus $\Delta$ itself is Fredholm acting between $x^{n/2-1} H^{2}_b(M) \to x^{n/2+1} L^2(M)$. Also $P$ is self-adjoint with respect to the measure induced by $g_b$, so its index is equal to zero, hence it is invertible if and only if its null space is trivial. This is therefore also true for $\Delta : x^{n/2-1} H^{2}_b(M) \to x^{n/2+1} L^2(M)$.It is also shown in \cite{tapsit}, section 5.25, that if $P f \in \CIdot(M)$, with $f$ in $L^2(M)$, then $f$ has an asymptotic expansion of the form
\begin{equation}
f \sim \sum_j x^{n-2+j} a_j \phi_j(\hat z), \quad
\label{f-exp}\end{equation}
where $\phi_j$ is a spherical harmonic with eigenvalue $j (j+n-2): \Delta_{S^{n-1}} \phi_j = j(j+n-2) \phi_j$. In particular, such a function tends to zero at infinity. It follows from this and from the maximum principle that there is no nontrivial solution to $\Delta f = 0$, with $f \in x^{n/2-1} H^{2}_b(M)$, because by \eqref{f-exp} $f$ would be a harmonic function tending to zero at infinity.  Hence we can solve $\Delta f = v$, $v \in \CIdot(M)$, where $f$ has an expansion \eqref{f-exp}. 

Let $f$ be as in the previous paragraph.
We first find a formal expansion for $u$ near the corner $\zf \cap \ff$ of $\Msp$. Coordinates near this corner are $x, y$ and $K = k/x$. Let us look for an expansion for $u$ of the form
\begin{equation}
u = \sum_j x^{n-2+j} \phi_j a_j(K), \quad  K = \frac{k}{x},
\end{equation}
where $a_j(0)$ is given by the expansion \eqref{f-exp} for $f$, so that $u | \zf = f$. 
The operator $\Delta + k^2$ may be written
\begin{equation}
(x^2 D_x)^2 + i(n-1) x^3 D_x + x^2  \Delta_{S^{n-1}} + k^2.
\end{equation}
Acting on the $j$th term this gives 
\begin{equation}
x^2 \Big( (xD_x)^2 + i(n-2) xD_x + \Delta_{S^{n-1}} + K^2 \Big).
\end{equation}
Here $D_x$ indicates the derivative keeping $k$ fixed. When we switch to using coordinates $(x, y, K)$,
 then we must replace $x D_x$ by $xD_x - K  D_K$, getting
\begin{equation}
x^2 \Big( (xD_x - K D_K)^2 + i(n-2) (xD_x - K D_K) + \Delta_{S^{n-1}} + K^2 \Big).
\end{equation}
Acting on the $j$th term, we may replace $\Delta_{S^{n-1}}$ by $j(j+n-2)$ and $xD_x$ by $-i(j+n-2)$, getting the operator
\begin{equation}
R_j \equiv x^2 \Big( ( K D_K)^2  -i(n-2+2j) (K D_K)  + K^2 \Big).
\end{equation}
The equation $R_j (a_j(K)) = 0$  has a smooth solution  for every $j$, with initial condition $a_j(0)$ determined by the coefficient in \eqref{f-exp}. We may cut this off with a cutoff function in $K$ whose derivative is supported where $K \in [1,2]$. The error term is then of the form
\begin{equation}
\sum_j x^{n+j} \phi_j \tilde b_j(K),
\end{equation} 
with $\tilde b_j$ supported in $[1,2]$. 

We now change to variables which are smooth at the other corner, $\bbdy \cap \ff$, namely $k$ and
 $\rho = x/k = 1/K$. The error term above may be written
\begin{equation}
\sum_j k^{n+j} \phi_j b_j(\rho),
\end{equation}
where $b_j$ is supported in $[1/2, 1]$. Let us try to solve it away with a series of the form
\begin{equation}
\sum_j k^{n-2+j} \phi_j c_j(\rho).
\end{equation}

Writing the operator in these new variables we get
\begin{equation}\begin{gathered}
(k \rho^2 D_\rho)^2 + i(n-1)k^2 \rho^3 D_\rho + x^2 \Delta_{S^{n-1}} + k^2 \\
= k^2 \Big( (\rho^2 D_\rho)^2 + i(n-1) \rho^3 D_\rho + \rho^2 \Delta_{S^{n-1}} + 1 \Big).
\end{gathered}\end{equation}
Let $c_j = e^{-1/\rho} e_j$. Then $e_j$ satisfies the equation 
\begin{equation}\begin{gathered}
  \Big((\rho^2 D_\rho -i)^2 + i(n-1) \rho (\rho^2 D_\rho - i) + \rho^2 j(j+n-2) + 1 \Big) e_j =  e^{1/\rho}b_j
\implies \\
 \Big( -2 \rho^2 \partial_\rho +(n-1)\rho +  \big( ( \rho^2 D_\rho)^2 + i(n-1) ( \rho^3 D_\rho) + 
 \rho^2 j(j+n-2) \big) \Big) e_j =  e^{1/\rho} b_j
\end{gathered}\end{equation}
This is a regular singular ODE with a solution of the form $e_j = \rho^{(n-1)/2} \tilde e_j$ where $\tilde e_j$ is smooth down to $\rho = 0$. This gives us a formal series in powers of $\rho_{\ff}$ at $\ff$ in which each term is uniformly rapidly decreasing with all derivatives at $\bbdy$ (i.e., as $\rho \to 0$). Borel summing at $\ff$, we get a formal solution that matches with $f$ to infinite order at $\ff$. Making a correction that vanishes to infinite order at $\zf$, in order to make $u$ agree exactly with $f$ at $\zf$, we get a function $u$ which satisfies all conditions  of the lemma.
\end{proof}

We now use this lemma to define a harmonic function on $M$ which will be key 
to the parametrix construction. We begin by choosing a smooth function
 $\phi : \RR^n \to \RR$ which is equal to $1$ for $|z|$ large and is supported in 
 $\{ |z| > 1 \}$. Using this we define functions $\phi_\pm$ on $M$, with 
 $\phi_\pm$ supported on the $\pm$ end of $M$, in the obvious way.
  Then let $u_\pm$ be the function given by Lemma~\ref{const} from the function 
  $v = -\phi_\pm$. It follows that
\begin{equation}
(\Delta + k^2) (e^{k z \cdot \omega} \phi_\pm + u_\pm) 
\in \rho_{\zf} \rho_{\bbdy}^\infty \rho_{\ff}^\infty \CI(\Mbar).
\label{bdy-1}\end{equation}
Moreover, 
\begin{equation}
\Phi_\pm = \phi_\pm + u_\pm | \zf \text{ is a harmonic function  equal to $1$ at $\partial \Mbar_\pm$ and $0$ at $\partial \Mbar_\mp$.}
\label{bdy-2}\end{equation}

We now define our parametrix. It is based on the resolvent kernel for $\RR^n$, but there is a crucial additional term ($G_3$ below) which corrects the leading order coefficient of the kernel at the face $\rbz$ (see the discussion of this coefficient in the Introduction). We now write $\phi_\pm$ to denote this function of the left variable o $M^2$ and $\phi_\pm'$ denote this function of the right variable on $M^2$.
Let $G_{\interior}(k)$ be a parametrix, modulo smoothing operators, for $(\Delta +k^2)^{-1}$ in the interior of $M$. We may assume that it is localized sufficiently close to the diagonal. We recall that the resolvent of the Laplacian on $\RR^n$ has the form \eqref{res}. Using this notation we define
\begin{equation}\begin{gathered}
\tilde G(k) =  G_1(k) + G_2(k) + G_3(k), \text{ where } \\
G_1(k) = 
\frac{e^{-k|z-z'|}}{|z-z'|^{n-2}} f_n(k |z-z'|) \Big( \phi_+ \phi_+' + \phi_- \phi_-' \Big) \\
G_2(k)  =  G_{\interior}(k) \Big( 1 - \phi_+ \phi_+' - \phi_- \phi_-' \Big) \\
G_3(k) = \frac{e^{-k |z'|}}{|z'|^{n-2}} f_n(k |z'|) \bigg(  u_+  \phi_+' 
+ u_- \phi_-' \bigg) .
\end{gathered}\label{G}\end{equation}

\section{Error term and resolvent}\label{error}
In this section, we correct the parametrix to the exact resolvent. The main point is that we obtain complete information about the regularity of the kernel of the error term on $\MMscsp$, and therefore of the resolvent itself on this space. This allows us to determine the regularity of $\Delta^{-1/2}$ and $d \Delta^{-1/2}$ and compute its behaviour to leading order at the boundary hypersurfaces of $\MMsc$.

Applying $\Delta + k^2$ to our parametrix (on the left), we get 
\begin{equation}
(\Delta + k^2) \tilde G(k) \equiv (\Delta + k^2) (G_1(k) + G_2(k) + G_3(k)) = \Id + \tilde E(k),
\end{equation}
where this equation defines $\tilde E(k)$. We may think of $\tilde E(k)$ either as a kernel on $\MMscsp$ or, by restricting to $M^2 \times [0, k_0]$ which is the interior of $\MMscsp$ together with the interior of $\zf$, a family of kernels parametrized by $k$ acting on functions on $M$.

By construction, the complete symbol of the diagonal singularity of $\tilde G(k)$ is the inverse, modulo symbols of order $-\infty$, of the complete symbol of $\Delta + k^2$.
Thus $\tilde E(k)$ is smooth at the diagonal. Also, we see due to the properties of $u_{\pm}$ that $\tilde E(k)$ vanishes to infinite order at $\bfc, \bfacez, \lb$, $\rb$ and $\lbz$.
 The crucial property of $\tilde E(k)$ is the order of vanishing at $\rbz$. 
 To calculate this, we need to determine the leading coefficient of the expansion of
  $\tilde G(k)$ at $\rbz$. These terms come from $G_1(k)$ and $G_3(k)$. 
  Since $|z-z'| = |z'| - z \cdot \hat z' + O(\rho_{\rbz}\rho_{\bfacez}^{-1})$ at $\rbz$, we have $k|z-z'| = k|z'| + O(\rho_{\rbz})$; note that  $k|z|$ vanishes on $\rbz$, while $k|z'|$ is finite in the interior of $\rbz$. Hence 
\begin{equation}
G_1(k) = \frac{e^{-k |z'|}}{|z'|^{n-2}} f_n(k |z'|)  \phi_\pm(z) + O(\rho_{\rbz}^{n-1}) \text{ at } \rbz.
\end{equation}
If we combine this with $G_3(k)$, then using \eqref{bdy-2} we see that the leading coefficient becomes
$$
\frac{e^{-k |z'|}}{|z'|^{n-2}} f_n(k |z'|) \Phi_\pm(z) + O(\rho_{\rbz}^{n-1})
 =  \rho_{\rbz}^{n-2} e^{-k |z'|} f_n(k|z'|) \Phi_\pm(z) + O(\rho_{\rbz}^{n-1})  \text{ at } \rbz.
$$
The leading term annihilated by the operator $\Delta + k^2$ (since $\Phi_\pm$ is harmonic and $k=0$ at $\rbz$), so the error term $\tilde E(k)$ is $O(\rho_{\rbz}^{n-1})$ at $\rbz$ --- an improvement of one order over what might be expected, and the main point of introducing the correction term $G_3(k)$. Thus, we have
\begin{equation}\label{Ek}
\tilde E(k) \in \rho_{\sf}^\infty \rho_{\bfc}^\infty \rho_{\bfacez}^\infty \rho_{\lb}^\infty \rho_{\lbz}^\infty \rho_{\rb}^{\infty} \rho_{\rbz}^{n-1} \CI(\MMscsp).
\end{equation}
Note that both $x = \ang{z}^{-1}$ and $x' = \ang{z'}^{-1}$ are smooth on $\MMscsp$, the former function vanishing simply at $\lb, \lbz, \bfc, \bfacez, \sf$ and the latter vanishing simply at $\rb, \rbz, \bfc, \bfacez, \sf$. Thus \eqref{Ek} implies that the kernel of $\tilde E(k)$ is $(\ang{z} \ang{z'})^{1-n}$ times a bounded function on $\MMscsp$. This implies that $\tilde E(k)$ is Hilbert-Schmidt, uniformly for $k \in [0, k_0]$, hence compact for each $k \in [0, k_0]$. Therefore $\Id - \tilde E(0)$ has finite dimensional null space and cokernel of the same dimension on $L^2(M)$.
We next show that we can modify our parametrix by the addition of a finite rank term so that the new error term is invertible for small $k$. The correction term will be
\begin{equation}
G_4(k) = \sum_{i=1}^N \phi_i \langle \psi_i, \cdot \rangle \, \quad N = \dim \null (\Id - \tilde E(0))
\label{g4-def}\end{equation}
where the $\phi_i$, $\psi_i$ are in $\CIdot(M)$ and independent of $k$. Since $\tilde E(0)$ maps into $\CIdot(M)$, the null space is contained in $\CIdot(M)$ and hence is independent of the choice of $l$. Thus we choose $\psi_i$ to span the null space of $\Id - \tilde E(0)$, and we would like to choose $\phi_i$ so that $\Delta \phi_i $ span a space supplementary to the range of $\Id + \tilde E(0)$. This is possible since $\Delta$ has trivial null space, and $\Delta$ is self-adjoint, hence the range of $\Delta$ on $\CIdot(M)$ is dense in $L^2(M)$. Choosing such $\phi_i$, we define $G_4(k)$ (which is actually independent of $k$) by \eqref{g4-def}. We now define
\begin{equation*}
G(k) = G_1(k) + G_2(k) + G_3(k) + G_4(k) = \tilde G(k) + G_4(k)
\end{equation*}
and define $E(k)$ by setting 
\begin{equation}
(\Delta + k^2) G(k) = \Id + E(k) \implies E(k) = \tilde E(k) + (\Delta + k^2) G_4(k);
\end{equation}
$E(k)$ enjoys all the properties of $\tilde E(k)$ listed above. In addition, since $E(0)$ is such that $\Id + E(0)$ is invertible, it follows that
actually $\Id + E(k)$ is invertible for all sufficiently small $k$; we assume that $k_0$ is chosen so that $\Id + E(k)$ is invertible for all $k \in [0, k_0]$.

We now analyze the inverse of $\Id + E(k)$. Let us write
$$
(\Id + E(k))^{-1} = \Id + S(k),
$$
where this equation defines $S(k)$. The decay of the kernel $E(k)$ at the boundary of $\MMscsp$ implies that $E(k)$ is Hilbert-Schmidt on $L^2(M)$. Hence $S(k)$ is also Hilbert-Schmidt. The regularity \eqref{Ek} of $\tilde E(k)$ on $\MMscsp$, and the fact that $x'/(x' + k) \in \CI(\MMscsp)$ vanishes simply 
at $\rb, \bfc, \sf$, imply that
$$E(k) \in x^N (x')^{n-1} \big( \frac{x'}{x' + k} \big)^N L^\infty(M^2 \times [0, k_0]) \text{ for all }N.$$
Using this and the formula $$S(k ) = E(k) + E(k)^2  + E(k) S(k) E(k)$$ shows that
\begin{equation}
S(k) \in x^N (x')^{n-1} \big( \frac{x'}{x' + k} \big)^N L^\infty(M^2 \times [0, k_0]) \text{ for all } N.
\label{Sk}\end{equation}

We are particularly interested in the kernel $G(k) S(k)$, which we shall call $G_5(k)$, since the addition of $G_5(k)$ will correct the parametrix $G(k)$ to the exact resolvent kernel.

\begin{lemma}\label{g5} Let $l = 0, 1, 2 \dots$. Then the kernel 
\begin{equation}
\int_0^{k_0} s_0(k) \nabla^{(l)} G(k) S(k) \, dk
\label{t5-int}\end{equation}
is in 
\begin{equation}\ang{z}^{-(n-1+l)} \ang{z'}^{-(n-1)} L^\infty(M^2) \cap \ang{z}^{-(n-2+l)} \ang{z'}^{-n} L^\infty(M^2).
\end{equation}
\end{lemma}

\begin{rem} Much more precise statements can be made about the kernels \eqref{Sk} and \eqref{t5-int}, for example by using Melrose's Pushforward Theorem \cite{RBMCalcCon}, which shows that these kernels are actually conormal, with respect to the boundary and the diagonal, on $\MMsc$. However, the $L^\infty$ statements will suffice for our purposes and are more straightforward to prove.
\end{rem}

\begin{proof}
Let us break up $G(k)$ into two parts $G(k) = \chi G(k) + (1 - \chi) G(k)$, where $\chi$ is as in Lemma~\ref{conormal}. Thus $\chi G(k)$ is a smooth family of scattering pseudodifferential operators,  while $(1 - \chi) G(k)$ has no singularity at the diagonal.

We first consider $(1 - \chi) G(k)$ which is localized away from the diagonal. Let $m_x$ denote the multiplication operator by $x=\ang{z}^{-1}$ on $M$. Then we have
$$
\nabla^{(l)} (1 - \chi) G(k)S(k) = \big( \nabla^{(l)} (1 - \chi) G(k)m_x^{n+l} \big) \big( m_x^{-(n+l)} S(k) \big).
$$
The kernel $\nabla^{(l)} G(k)$ decays to order $n-2+l$ at $\lbz$ and to order $\infty$ at $\lb$ and $\bfc$. If we multiply this kernel by $\ang{z'}^{-(n+l)}$, which corresponds to composing with $m^{n+l}$ on the right, then it also decays to order $n+l$ at  $\bfacez$ and $\sf$. This means that we can write
\begin{equation}
\nabla^{(l)} (1 - \chi) G(k)m_x^{n+l} \in x^{n-2+l} \big( \frac{x}{x+k} \big)^2 L^\infty(M^2 \times [0, k_0]), \quad x = \ang{z}^{-1}
\end{equation}
since $x$ is a product of boundary defining functions for $\lb, \lbz, \bfc, \bfacez, \sf$, and $x/x+k$ vanishes to first order at $\lb$ and $\bfc$. In a similar way, using \eqref{Sk}, we find that
\begin{equation}
m_x^{-(n+l)} S(k) \in (x')^{n-1} \big( \frac{x'}{x'+k} \big)^2 L^1(M; L^\infty(M \times [0, k_0]));
\end{equation}
note that composing with $m_x^{-(n+l)}$ on the left is harmless here because the kernel $S(k)$ vanishes to infinite order on every boundary hypersurface where $x^{-(n+1)}$ blows up. For this same reason the kernel is $L^1$ in the left variable, uniformly in the right variable and in $k$. It follows that the composition
\begin{equation}
\nabla^{(l)} (1 - \chi) G(k)S(k) \in x^{n-2+l} (x')^{n-1} \big( \frac{x}{x+k} \big)^2 \big( \frac{x'}{x'+k} \big)^2 L^\infty(M^2 \times [0, k_0]).
\label{Linfty}\end{equation}
Now we integrate in $k$. If we ignore the $(x/x+k)^2$ factor (which is bounded) then we find that
$$
\int_0^{k_0} s_0(k) \nabla^{(l)} (1 - \chi) G(k)S(k) \, dk \in x^{n-2+l} (x')^n L^\infty(M^2),
$$
because
\begin{equation*}
\int_0^{k_0} \big( \frac{x'}{x'+k} \big)^2 \, dk \leq 
\int_0^{\infty} \big( \frac{x'}{x'+k} \big)^2 \, dk = C x', \
C = \int_0^\infty \big( \frac{1}{1+\overline{k}} \big)^2 \, d\overline{k}.
\end{equation*}
In exactly the same way we show that 
$$
\int_0^{k_0} s_0(k) \nabla^{(l)} (1 - \chi) G(k)S(k) \, dk \in x^{n-1+l} (x')^{n-1} L^\infty(M^2).
$$

Finally we consider the integral \eqref{t5-int} with $G(k)$ replaced by $\chi G(k)$. We may regard $\nabla^{(l)} \chi G(k)$ as a smooth family of scattering pseudodifferential operators, and $S(k)$ as an element of
\begin{equation}
\rho_{\rb}^{\infty} \rho_{\rbz}^{n-1} \CIdot(M; \CI(\Msp)).
\label{S-space}\end{equation}
Since scattering pseudodifferential operators map $\CIdot(M)$ to itself continuously, it follows that $\nabla^{(l)} \chi G(k) S(k)$ is also an element of the space \eqref{S-space}. Performing the $k$ integral we get an extra vanishing factor at $\rb \subset \MMsc$, yielding $x^\infty (x')^n L^\infty(M^2)$, which proves the Lemma for this piece. This completes the proof.
\end{proof}

\section{Riesz Transform}\label{RT}
Recall that in Section~\ref{para} we split $\Delta^{-1/2} = g_0(\Delta) + g_1(\Delta)$, where $d g_1(\Delta)$ was bounded from $L^p$ to $L^p$ for all $1<p<\infty$, and
\begin{equation}
g_0(\Delta) = \frac{2}{\pi} \int_0^\infty s_0(k) (\Delta + k^2)^{-1} \, dk.
\end{equation}
Hence, it remains to analyze $d g_0(\Delta)$.
Let us decompose
\begin{equation}
(\Delta + k^2)^{-1} = G_1(k) + G_2(k) + G_3(k) + G_4(k) + G_5(k)
\end{equation}
as in the previous section and write $d g_0(\Delta) = T_1 + T_2 + T_3 + T_4 + T_5$ correspondingly.

The easiest kernel to deal with is $T_4$; this kernel is in $\CIdot(M^2)$, hence is bounded from $L^p$ to $L^p$ for $1 \leq p \leq \infty$. The kernel $T_2$ is bounded on $L^p$ for $1 < p < \infty$ because it is a classical zero order pseudodifferential operator with proper support; see Chapter VI, section 5 of \cite{Stein}. The kernel $T_1$ we decompose further as as $T_1 = T_{1,1} + T_{1,2} + T_{1,3}$, where \begin{equation}\begin{gathered}T_{1,1} = \Big( d_z \frac{e^{-k|z-z'|}}{|z-z'|^{n-2}} f_n(k |z-z'|) \Big) \Big( \phi_+ \phi_+' + \phi_- \phi_-' \Big) , \\T_{1,2} = \chi \Big( \frac{e^{-k|z-z'|}}{|z-z'|^{n-2}} f_n(k |z-z'|) \Big) \Big(  (d\phi_+) \phi_+' +  (d\phi_-) \phi_-' \Big), \\T_{1,3} = (1 - \chi)  \Big( \frac{e^{-k|z-z'|}}{|z-z'|^{n-2}} f_n(k |z-z'|) \Big) \Big(  (d\phi_+) \phi_+' +  (d\phi_-) \phi_-' \Big)\end{gathered}\end{equation}and $\chi$ is as defined above Lemma~\ref{conormal}. It is clear that $T_{1,1}$ is bounded on $L^p$ for $1 < p < \infty$, because t
 he Riesz kernel on $\RR^n$ has this property. Also, $T_{1,2}$, like $T_2$, it is a classical zero order pseudodifferential operator with proper support, hence bounded on all $L^p$.

We next consider $T_5$.  Lemma~\ref{g5}, with $l=1$, shows that $T_5$ is
in $L^p(M; L^{p'}(M))$, where $p'^{-1} = 1 - p^{-1}$, for all $p \in (1, \infty)$, which implies that $T_5$ is bounded on $L^p$ for $1 < p < \infty$. Thus we are left with $T_{1,3} + T_3$.

\begin{lemma} The kernel of $T_{1,3} + T_3$ such that 
\begin{equation}
T_{1,3} + T_3 \in \rho_{\rb}^{n-1}  \rho_{\lb}^n \big( \rho_{\sf} \rho_{\bfc}  \big)^{2n-2} \CI(\MMsc).
\label{T}\end{equation}
Moreover, the leading coefficient of $T_{1,3} + T_3$ at $\rb$ is a constant times $d \Phi_\pm$.
\end{lemma}

\begin{proof}
Let us first consider the kernel of $T_{1,3}$ near $\rb$ and away from $\bfc$. This given by
\begin{equation*}
\int_0^{k_0} \bigg(
\frac{e^{-k|z-z'|}}{|z-z'|^{n-2}} f_n(k |z-z'|) \Big( d \phi_+ \phi_+' + d \phi_- \phi_-' \Big)
\bigg) dk;
\end{equation*}
note that $T_{1,3}$ is supported away from $\lb$, $\bfc$ and $\sf$ since the support of $d \phi$ is compact. It is a smooth function of $z$ and $x', y', k/x'$ which is rapidly decreasing in $k/x'$. It vanishes to order $(x')^{n-2}$ at $x' = 0$; note that $x'$ is a boundary defining function for $\rb$ in this region. Moreover, it is given by
$$
(x')^{n-2} e^{-k/x'} f_n(k/x') \Big( d \phi_+ \phi_+' + d \phi_- \phi_-' \Big) + O((x')^{n-1}).
$$
Changing variable of integration to $k/x'$ and taking into account $dk = x' d(k/x')$ we see that the integral is
$$
C_n (x')^{n-1} \Big( d \phi_+ \phi_+' + d \phi_- \phi_-' \Big) + O((x')^{n})
$$
at $\rb$. If we do the analogous calculation for $T_3$ and add the results we find that the kernel of $T_{1,3} + T_3$ is given by
$$
C_n (x')^{n-1} \Big( d \Phi_+ \phi_+' + d \Phi_- \phi_-' \Big) + O((x')^{n})
$$
at $\rb$. This proves the last statement of the lemma.

A similar computation can be done for $T_3$ at $\lb$, but now the result vanishes to order $n$ at the left boundary, because the derivative $d$, which is applied to the left variable of the kernel, increases the order of vanishing by $1$ at the left boundary.

Consider next the kernel $T_3$ near the triple intersection of $\rb$, $\rbz$ and $\bfacez$. 
In this case, local boundary defining functions are $\rho_{\rb} = x'/k$, $\rho_{\bfacez} = x$ 
and $\rho_{\rbz} = k/x$.
 We claim that the kernel is actually a smooth function of $x'/k$, $x'/x$, $x$, $y$ 
 and $y'$ in this region, which is a stronger statement, since $x'/x = \rho_{\rb} \cdot \rho_{\rbz}$. 
 To see this, note that the kernel of $T_3$ is equal to $e^{-1/\rho_{\rb}}$ times a $\CI$ function 
 on $\MMscsp$. Generally, if $h(u,v)$ is any smooth function of $u$ and $v$, $u, v \geq 0$, then $e^{-1/u} h(u,v)$ is a smooth function of $u$ and $v/u$. In other words, the function $\tilde h(u, w) = e^{-1/u} h(u, w/u)$ is smooth. (This is easily checked directly by differentiating $\tilde h$; inverse powers of $u$ are harmless due to the $e^{-1/u}$ factor.) Now let $u = \rho_{\rb}$ and $v = \rho_{\rbz}$, and treat the other coordinates as parameters, and the claim follows. 

The kernel $T_3$ vanishes to order $n-2$ at $\rbz$ and $2n-3$ at $\bfacez$. We change variable of integration to $k/x'$ as before, and the change of measure $dk = x'd(k/x')$ gives us additional vanishing at \emph{both} $\rbz$ and $\bfacez$, since $x' = (x'/x) x$ vanishes at both $\rbz$ and $\bfacez$. Thus the result is a smooth function of $(x'/x, x, y, y')$ which vanishes to order $n-1$ at $\rb = \{x'/x = 0\}$ and order $2n-2$ at $\bfc = \{x = 0\}$, which verifies the statement of the lemma near the corner $\rb \cap \bfc \subset \MMsc$. The other regions of $\MMsc$ are treated similarly.
 \end{proof}

This lemma implies that, for $1 < p < n$, $T_{1,3} + T_3$ is an element of  $L^p(M; L^{p'}(M))$. Moreover, for $p \geq n$, this is not true since the function $(x')^{n-1}$ is not in $L^{p'}$ then, and the coefficient of $(x')^{n-1}$ is $d \Phi_\pm$ which does not vanish identically. Therefore $T_{1,3} + T_3$ cannot be applied to any bounded function equal to $x(\log x)^{-1}$ near infinity, which lies in $L^p$ for $p \geq n$. This completes the proof of Theorem~\ref{main}.

\begin{rem} If $M$ has one Euclidean end then the same argument shows that the Riesz transform is bounded on $L^p$ for all $1 < p < \infty$. In this case, the parametrix $\tilde G(k)$ can be taken to be (compare with \eqref{G})
$$
\frac{e^{-k|z-z'|}}{|z-z'|^{n-2}} f_n(k |z-z'|)  \phi \phi' +
G_{\interior}(k) \big( 1 - \phi \phi' \big) +
\frac{e^{-k |z'|}}{|z'|^{n-2}} f_n(k |z'|)   (1-\phi)  \phi' .
$$
In this case the role of $\Phi_\pm$ in the computation above is played by the constant function $1$. The argument is the same as above, except that the gradient of $1$ vanishes so that we get $\rho_{\rb}^n$ instead of $\rho_{\rb}^{n-1}$ in \eqref{T} (as outlined in the introduction), leading to the boundedness for all $p$ strictly between $1$ and $\infty$. 
\end{rem}


\section{Heat kernel}\label{heat}
As part of the analysis of the heat kernel we analyzed the structure of the resolvent $(\Delta + k^2)^{-1}$ for real $k$, including an analysis of the asymptotics of its kernel when $k \to 0$. This analysis remains valid for any cone $\{ k = i\lambda \mid \Imag \lambda \geq \epsilon \Real \lambda \}$ for any $\epsilon > 0$. We can use this  to obtain information about the heat kernel $H(t, z, z')$ of $e^{-t \Delta}$ on $M$ via the contour integral
\begin{equation}
e^{-t \Delta} = \frac1{2\pi i} \int_\Gamma e^{-t \lambda^2} \big( \Delta - \lambda^2 \big)^{-1} 2\lambda \, d\lambda
\label{contour-int}\end{equation}
where $\Gamma$ is the contour $\{ \lambda = s e^{-i\pi/12} \cup 
\lambda = s e^{i\pi/12} \mid s \in \RR^+ \}$. 

Let us focus on the heat kernel in the following asymptotic regime: We fix a point $z \in M$, which we think of as being in the `compact part' of $M$ (where the metric is not flat), and fix an end of $M$ and a point $\omega \in S^{n-1}$ which we think of as a point at infinity for this end. Consider the behaviour of the heat kernel $H(t,z, z')$ where $z' = r' \omega$ and $t \to \infty$, $r' \to \infty$ so that $\sqrt{t}/r'$ approaches a finite positive limit $\sigma$.

\begin{proposition} Assume that $M$ has Euclidean ends, with the number of ends at least two.
Under the limiting process described above, 
$t^{n/2} \nabla_{z}^{(l)} H(t, z, z')$ approaches a limit, for any value of $l$. Indeed
\begin{equation}
\lim_{t \to \infty} t^{n/2} \nabla_{z}^{(l)} H(t, z, z') = (4\pi)^{-n/2} \ e^{-1/4\sigma^2} \nabla_{z}^{(l)} \Phi(z),\  \sigma  = \frac{\sqrt{t}}{r'} > 0 \text{ fixed.}
\label{heat-limit}\end{equation}
where $\Phi$ is the harmonic function which tends to $1$ at the given end and tends to $0$ at all other ends. In particular, we have a lower bound on the derivatives of the heat kernel for large time:
\begin{equation}
\sup_{z, z' \in M} \big|  \nabla_z^{(l)} H(t, z, z') \big| \geq c_l t^{-n/2}, \, \text{ for some } \, c_l > 0, \quad t \geq 1. 
\end{equation}
\end{proposition}

\begin{rem} For $l=0$ this result is not surprising. The point of this proposition is that \emph{taking derivatives in the $z$ variable gives no additional decay in the heat kernel} (in this asymptotic regime). This contrasts with Euclidean space where each additional derivative gives additional decay of $t^{-1/2} = (\sigma r')^{-1}$.
\end{rem}

\begin{proof} The $k$th $z$-derivative of the heat kernel is given by the contour integral \eqref{contour-int} 
with the resolvent replaced by the $k$th $z$-derivative of the resolvent. Clearly, to prove the theorem we only have to consider the kernel of the resolvent in a neighbourhood of $\rb$ and $\rbz$.

Near the interior of $\rbz$, and away from $\rb$ the function $\Lambda = \lambda/x'$ is a smooth function, which goes to infinity at $\rb$; in fact, $\Lambda^{-1}$ is a boundary defining function for $\rb$. In this integral \eqref{contour-int}, the term $e^{-t\lambda^2} = e^{-\sigma^2 \Lambda^2}$ therefore vanishes together with all its derivatives at $\rb$, since $\sigma > 0$ by assumption, which means that we may ignore the expansion of the resolvent at $\rb$. Hence to find the asymptotics the heat kernel near $\rb$ in this regime we only need to consider the expansion of the resolvent at $\rbz$ (up to a correction that vanishes to infinite order as $r' \to \infty$).

Using the $L^\infty$ bounds \eqref{Linfty},
 we may write the $k$th derivative of the resolvent kernel in the form
\begin{equation}\begin{gathered}
K_0(z, y', \Lambda) + K_1(z, y', \Lambda, r'), \\
K_0 = (r')^{-(n-2)} e^{i \Lambda} f_n(\Lambda) \nabla_{z}^{(l)} \Phi(z), \quad  K_1 = O((r')^{-(n-1)})
\end{gathered}\end{equation}
in the region of interest. 
Let us first substitute $K_0$ for the resolvent into the integral \eqref{contour-int}.  Thus we want to compute the limit
\begin{equation}
\lim_{t \to \infty} t^{n/2} \frac1{\pi i} \int_\Gamma e^{-t \lambda^2} (r')^{-n+2} e^{i \Lambda} f_n(\Lambda) \nabla_{z}^{(l)} \Phi(z) \lambda \, d\lambda.
\label{contour-int-2}\end{equation}
Substituting $\lambda = (r')^{-1} \Lambda$ and $t = \sigma^2 (r')^2$, and using $\lambda d\lambda = (r')^{-2} \Lambda d\Lambda$, we get
\begin{equation}
\lim_{t \to \infty} \frac1{\pi i}  \sigma^{n}   \nabla_{z}^{(l)} \Phi(z)  \int_\Gamma e^{-\sigma^2 \Lambda^2}  e^{i \Lambda} f_n(\Lambda) \Lambda \, d\Lambda.
\label{contour-int-3}\end{equation}
Taking the limit is trivial, since \eqref{contour-int-3} is independent of $t$.
To perform the integral, consider the case of $\RR^n$, with kernel $(\Delta - \lambda^2)^{-1}(z,z')$ with $z$ fixed to be the origin. This gives rise to an integral
\begin{equation}
\frac1{\pi i}  (r')^{-n} \int_\Gamma e^{- \sigma^2 \Lambda^2}  e^{i \Lambda} f_n(\Lambda) \Lambda \, d\Lambda
\label{contour-int-4}\end{equation}
which is equal to 
$$
(4 \pi t)^{-n/2} e^{-(r')^2/4t} = (4 \pi t)^{-n/2} e^{-1/(4\sigma^2)}.
$$
Multiplying through by $(r')^n$ gives 
\begin{equation*}
\frac1{\pi i}   \int_\Gamma e^{- \sigma^2 \Lambda^2}  e^{i \Lambda} f_n(\Lambda) \Lambda \, d\Lambda = (4\pi)^{-n/2} \sigma^{-n} e^{-1/(4\sigma^2)}.
\end{equation*}
Hence, \eqref{contour-int-3} is equal to 
\begin{equation}
(4\pi)^{-n/2}  e^{-1/(4\sigma^2)}    \nabla_{z}^{(l)} \Phi(z) ,
\end{equation}
which is the right hand side of \eqref{heat-limit}. 
If we now substitute $K_1$ for the resolvent in \eqref{contour-int}, which vanishes to an additional order as $r' \to \infty$ as compared to $K_0$, then the integral also vanishes to an additional order, giving a zero contribution to the limit \eqref{heat-limit}. This proves the proposition.
\end{proof}
It is also of interest to compute the leading behaviour of the heat kernel $H(t, z, z')$ as $t \to \infty$ and as $z, z'$ both tend to infinity, but along different ends. Suppose that $z = r \omega$, where $\omega \in S^{n-1}_-$ is fixed and that $z' = r' \omega'$, $\omega' \in S^{n-1}_+$ is fixed, and suppose further that $\sqrt{t}/r \to \sigma, \sqrt{t}/r' \to \sigma'$ where $\sigma, \sigma' \in (0, \infty)$. \begin{proposition} Under this asymptotic regime, the limit \begin{equation}\lim_{t \to \infty} t^{n-1} H(t, z, z') = q(\sigma, \sigma'),  \quad \frac{\sqrt{t}}{r} \to \sigma, \ \frac{\sqrt{t}}{r'} \to \sigma'\end{equation}exists and is finite. Hence in this asymptotic regime the heat kernel has $t^{-n+1}$ decay as $t \to \infty$.\end{proposition}\begin{rem} For $n \geq 3$ this is faster than the usual $t^{-n/2}$ decay. Hence Gaussian lower bounds do not hold for the heat kernel on $M$. This was observed in \cite{BCF}, and can be heuristically explained in terms of Bro
 wnian motion on $M$. Here we give an explicit quantitative description of the failure of this lower bound. \end{rem}\begin{proof}We shall perform a similar computation as in the proof of the previous proposition. Since $|z|/|z'| \to \sigma'/\sigma$ under this limiting regime, and $z, z'$ go to infinity along different ends, we end up at the `anti-diagonal' part of $\bfc$. Hence we need to consider the resolvent kernel near the anti-diagonal part of $\bfc$ and $\bfacez$, where $y \in S^{n-1}_-$ and $y' \in S^{n-1}_+$. It is the $G_3(k)$ term which is important here; we need the leading behaviour of $u_+$ at the negative end. It is not hard to show that $$u_+ = A |z|^{-n+2} e^{-k|z|} f_n(k|z|) + O(|z|^{-n+1}) \text{ for some } A > 0. $$at this end. Indeed, the harmonic function $\Phi_+$ is equal to $A' |z|^{-n+2} + O(|z|^{-n+1})$ as $z \to \infty$ along this end, for some $A' > 0$. The leading coefficient $a_0(K)$ from Lemma~\ref{const} must then be equal to a constant times $
 e^{-K} f_n(K)$, which follows readily from the fact that $|z|^{-n+2} e^{-k|z|} f_n(k|z|)$ satisfies the equation $(\Delta + k^2) u = 0$. The specific structure of the parametrix $G(k)$, together with the estimate \eqref{conorm} with $l=0$, shows that in this region the resolvent kernel may be written as a sum\begin{equation}\begin{gathered}
K_0(y, y', r', \Lambda) + K_1(y, y', r', \Lambda), \\
K_0 = A r^{-n+2} (r')^{-n+2} e^{i \lambda r} f_n(\lambda r)e^{i \lambda r'} f_n(\lambda r') , \quad  K_1 = O((r')^{-2(n-1)})
\end{gathered}\end{equation}
in the region of interest. Substituting $K_0$ for the resolvent into the integral \eqref{contour-int}, we obtain\begin{equation}
\lim_{t \to \infty} t^{n-1} \frac1{\pi i} \int_\Gamma e^{-t \lambda^2} A r^{-n+2} (r')^{-n+2} e^{i \lambda r} f_n(\lambda r)e^{i \lambda r'} f_n(\lambda r') \lambda \, d\lambda.
\label{contour-int-22}\end{equation}Let $\alpha = \sigma'/\sigma = \lim r/r'$.
Substituting $\lambda = (r')^{-1} \Lambda$ and $t = (\sigma')^2 (r')^2$, and using $\lambda d\lambda = (r')^{-2} \Lambda d\Lambda$, we get
\begin{equation}
\lim_{t \to \infty}  t^{n-1} (r')^{-2(n-1)} \alpha^{-n+2}   \frac1{\pi i}   \int_\Gamma e^{-(\sigma')^2 \Lambda^2}  e^{i \Lambda} f_n(\Lambda) e^{i \alpha \Lambda} f_n(\alpha \Lambda) \Lambda \, d\Lambda = C(\alpha, \sigma').
\label{contour-int-33}\end{equation}Thus the limit exists and is finite, when $K_0$ is substituted for the resolvent. As in the previous proof, when $K_1$ is substituted for the resolvent the limit is zero, since $K_1$ decays to an additional order at infinity. This completes the proof. \end{proof}

\section{Riesz transform and $L^p$ cohomology}\label{coh}Here  $(M^n,g)$ is a complete Riemannian manifold of dimension $n$.

  We want here to discuss some consequence of the boundedness of the Riesz transform on $L^p$ for some $p> 2$  for the $L^p$ cohomology. On $(M,g)$, the space of $L^2$ differential forms admits the Hodge decomposition 
  $$L^2(T^*M)=\cH^1(M)\oplus \overline{d C^\infty_0(M)}\oplus \overline{d^*
C^\infty_0(\Lambda^2 T^*M)},$$
  where $\cH^1(M)=\{\alpha\in L^2(T^*M), d \alpha=0=d^*\alpha\}$
(see \cite{DR}).
  Let us recall now the definition of reduced $L^p$-cohomology:
for $p\ge 1$ , the first space of reduced $L^p$ cohomology of $(M,g)$ is
  $$H^1_p(M)=\frac{\{\alpha\in L^p(T^*M),
  d\alpha=0\}}{\overline{dC^\infty_0(M)}} \ ,$$
  where we take the closure in $L^p$. 
The first space of reduced $L^2$ cohomology can be identified with $\cH^1(M)$.
 As noticed in \cite{ACDH}, if we assume that for  some $p\ge 2$ the Riesz
 transform $T:=d\Delta^{-1/2}$ is  bounded on $L^p$ and on $L^{p/(p-1)}$, then the Hodge
 projector $$P=d\Delta^{-1}d^*=TT^*:L^p(M; T^*M)\cap L^2(M; T^*M)\longrightarrow
 L^2(M; T^*M)$$ i.e. the orthogonal projector of $L^2(M; T^*M)$ onto the space of
 `exact forms'  extends by continuity to a bounded operator  $$P: L^p(M; T^*M)\longrightarrow L^p(M; T^*M).$$

We assume now that $(M^n,g)$ is a complete Riemannian manifold, $n \geq 2$, satisfying the Nash inequality
\begin{equation}
\label{nash}\mu \left(\int_M f^2 d\vol\right)^{1+2/n}\le \left(\int_M |f|
d\vol\right)^{4/n} \int_M |df|^2 d\vol,
\end{equation}
for all $f\in C^\infty_0(M)$ and some $\mu>0$, and that the volume growth of geodesic balls is uniformly bounded:
\begin{equation}
\label{volumgro}
\forall x\in M,\, \forall r>0,\ \vol B(x,r)\le C r^n. 
\end{equation}
It follows from \cite{Ca} that (\ref{nash}) implies a matching lower bound:
\begin{equation}
\label{volumlo}
\forall x\in M,\, \forall r>0,\ \vol B(x,r)\ge c r^n.
\end{equation}
Note that (\ref{nash}) easily implies the Faber-Krahn inequality :
\begin{equation}\label{fk}
\ \lambda_1(\Omega) \ge \mu
\left( \vol \Omega\right)^{-2/n},
\end{equation}
for all  $\Omega\subset M$ with finite measure,
where $$\lambda_1(\Omega)=\inf\left\{\frac{ \int_\Omega |df|^2d\vol}{\int_\Omega f^2d\vol} , f\in
C^\infty_0(\Omega)\setminus\{0\}\right \}$$
is the first eigenvalue for the Laplacian on $\Omega$ for the Dirichlet
boundary conditions (in fact, (\ref{nash}) and (\ref{fk}) are equivalent, see \cite{G}).
Also, if $n>2$,  (\ref{nash}) is equivalent  
 to  the Sobolev inequality :
 \begin{equation}
 \label{sobo}
 \nu \left(\int_M |f|^{\frac{2n}{n-2}} d\vol\right)^{1-2/n}\le  \int_M
 |df|^2 d\vol,\ \forall f\in C^\infty_0(M),
\end{equation}
for some $\nu>0$
(see for instance \cite{BCLS}).

According to \cite{CD}, we know that on $(M,g)$ the
Riesz transform  is
bounded on $L^q$ for $q\in]1,2]$. 
Hence if we assume that for  some $p\ge 2$ the Riesz transform is also bounded on $L^p$ then,
according to the above remark,
the Hodge projector 
$$P=:L^p(M; T^*M)\cap L^2(M; T^*M)\longrightarrow L^2(M; T^*M)$$
extends by continuity to a bounded
operator  $$P: L^p(M; T^*M)\longrightarrow L^p(M; T^*M).$$
  \begin{lemma}\label{last} Under the hypotheses
  (\ref{nash}, \ref{volumgro}), if the Riesz  transform is bounded in $L^p$ for some $p>2$, then
   $P\left(L^p(M; T^*M)\right)$ is the closure in $L^p(M; T^*M)$ of
   $d C^\infty_0(M)$.
 \end{lemma} \proof

According to \cite{CKS}, the Nash inequality implies that the semigroup
$e^{-t\Delta}$ satisfies the bound 
\begin{equation} 
\label{heatmapping} 
\|e^{-t\Delta}\|_{L^1\to L^\infty}\le Ct^{-n/2}, \forall\,t>0. 
\end{equation} 
A result of N. Varopoulos (\cite{V}) then implies the following mapping property for $q\in ]1,n[$:
 \begin{equation} 
\label{mapping} \Delta^{-1/2}:L^q(M)\rightarrow
L^{qn/(n-q)}(M).\end{equation}
  In order to prove the lemma, we have to show that if $\alpha\in C^\infty_0(M; T^*M)$ then $P\alpha$ can be approximated
   in $L^p$ by a sequence of elements  of $d C^\infty_0(M)$. We seek a sequence $\chi_k$
    of smooth functions with compact support such that
    $$L^p \operatorname{-}\lim_{k\to\infty} d(\chi_k\Delta^{-1} d^*\alpha)=P\alpha.$$
    Since we assume that the Riesz transform is bounded in $L^p$, we know that its
    adjoint $$\Delta^{-1/2}d^*:L^{p/(p-1)}(M; T^*M)\rightarrow L^{p/(p-1)}(M)$$ is
    bounded.  Hence we have $\Delta^{-1/2}d^*\alpha\in L^{p/(p-1)}(M)$.  Note that the condition
$p/(p-1)<n$ is satisfied since we are assuming that $n\ge 2$, in which case  $p> 2\ge  n/(n-1)$. Thus, by 
\eqref{mapping},
$$\Delta^{-1/2}\Delta^{-1/2}d^*\alpha=\Delta^{-1}d^*\alpha\in
L^{\frac{pn}{n(p-1)-p}}(M).$$Choose a point $o \in M$ and choose
 $$\chi_k(x)=
 \left\{\begin{array}{lll}
 1& {\rm if} & x\in B(o,k)\\
0&{\rm if} & x\not \in B(o,2k)\\
\end{array}\right.$$
$${\rm with}\  \|d\chi_k\|_{L^\infty}\le C/k.$$
Since $$d(\chi_k\Delta^{-1} d^*\alpha)=\chi_k(d\Delta^{-1}
d^*\alpha)+d\chi_k(\Delta^{-1} d^*\alpha )$$
and 
$$\chi_k(d\Delta^{-1}
d^*\alpha)=\chi_kP\alpha$$
obviously tends to $P\alpha$ in $L^p$ as $k\to\infty$, we need only to show
that $$\lim_{k\to\infty} \|d\chi_k(\Delta^{-1} d^*\alpha) \|_{L^p}=0.$$
We know that $\varphi=\Delta^{-1} d^*\alpha $ is harmonic outside a big ball
$B(o,R_0)$ containing the support of $\alpha$. Now the  Faber-Krahn
inequality (\ref{fk}) implies a mean value inequality for harmonic functions (see  \cite{G1}, Lemma 6.9) which yields:
$$|\varphi(x)|\le
C(\mu,n,p)r^{-\frac{(p-1)n-p}{p}}\|\varphi\|_{L^{\frac{pn}{n(p-1)-p}}(B(x,r))},$$
provided that $B(x,r)\subset M\setminus B(o,R_0)$ . 
In particular if $\rho(x)=\dist(x,o)-R_0>0$ we obtain $$|\varphi(x)|^p \le \frac{C}{ \rho(x)^{(p-1)n-p}}.$$
Hence we finally obtain, if say $k\ge 2R_0$:
$$ \|d\chi_k(\Delta^{-1} d^*\alpha) \|_{L^p}^p\le C\frac{\vol
B(o,2k) }{k^{(p-1)n}} \le C k^{(2-p)n}$$
which indeed goes to zero when $k\to\infty$.
We have proved that
$$PL^p(M; T^*M))\subset \overline{dC^\infty_0(M)}.$$
The converse inclusion follows from the fact that
$$dC^\infty_0(M)\subset PL^p(M; T^*M))$$
and that $P$, being a bounded projector, has a closed range.
\endproof

As a consequence of Lemma \ref{last},  if the assumptions (\ref{nash},\ref{volumgro}) are satisfied and if  the Riesz transform is bounded on $L^p$ for some $p>2$,
then $H^1_p(M)$ can be identified with:$$ \{\alpha \in L^p(M; T^*M) \mid
d\alpha=0\ {\rm  and}\  P\alpha=0\}.$$ Moreover we also have
\begin{lemma}  Under the hypotheses
  (\ref{nash}, \ref{volumgro}), if the Riesz  transform is bounded in $L^p$ for some $p>2$, then
$$\{\alpha\in L^p(M; T^*M) \mid d^*\alpha=0\}=\{\alpha\in L^p(M; T^*M) \mid P\alpha=0\}.$$
\label{coclosed}
\end{lemma}
\proof As a matter of fact, we have 
$\{\alpha\in L^p(M; T^*M), P\alpha=0\}=\Imag (\Id -P)$. The density of $C_0^\infty(M; T^*M)$ in $L^p(M; T^*M)$ and the boundedness of $P$ in $L^p$ imply that $(\Id
-P)(C_0^\infty(M; T^*M))$ is dense in $\Imag (\Id -P)$.
But we have $$(\Id -P)(C_0^\infty(M; T^*M))\subset \{\alpha\in L^p(M; T^*M) \mid d^*\alpha=0\}.$$ The latter space is closed,  hence we have the inclusion $$\{\alpha\in L^p(M; T^*M) \mid P\alpha=0\}\subset \{\alpha\in L^p(M; T^*M) \mid d^*\alpha=0\}.$$
Now assume that $\alpha\in L^p(M; T^*M)$ is coclosed. We define a sequence of  cutoff functions by
\begin{equation}\chi_k(x)=
     \left\{\begin{array}{lll}
     1& {\rm if} & x\in B(o,k)\\
     \frac{ \log\left(k^2/\dist(x,o)\right)}{\log k} & {\rm if} & x\in B(o,k^2)\setminus B(o,k)\\
     0&{\rm if} & x\not \in B(o,k^2)\\
     \end{array}\right.\label{cutoffk}\end{equation}
Then $L^p-\lim_{k\to\infty}\chi_k\alpha=\alpha$ but now $\chi_k \alpha\in
L^2$ and $P\chi_k\alpha=-T\Delta^{-1/2}({\rm int}_{\nabla\chi_k}\alpha)$,
because $d^*(\chi_k \alpha)=-{\rm int}_{\nabla\chi_k}\alpha$ (here ${\rm int}_{\nabla\chi_k}\alpha$ denotes the contraction of $\alpha$ with the vector field $\nabla \chi_k$).
Take $l\in ]1,n[ $ with $$\frac1l=\frac1p+\frac1n.$$ 
Then we have 
$$\left\|{\rm int}_{\nabla\chi_k}\alpha\right\|_{L^l}\le
\left\|\alpha\right\|_{L^p(M\setminus B(o,k))}\, \left\|\nabla\chi_k\right\|_{L^n}.$$
But an easy computation leads to 
$$\left\|\nabla\chi_k\right\|^n_{L^n}\le C (\log k)^{-(n-1)}.$$ With  \eqref{mapping}, we have
$$\lim_{k\to \infty} 
\left\|\Delta^{-1/2}({\rm int}_{\nabla\chi_k}\alpha)\right\|_{L^p}=0$$
hence by continuity of $T$ we obtain $P\alpha=0$.
\endproof
In particular we obtain a Hodge-de Rham interpretation of the $L^p$ cohomology:
\begin{proposition}\label{HdR} Under the hypotheses
  (\ref{nash}, \ref{volumgro}), if the Riesz  transform is bounded in $L^p$ for some $p>2$, then
  $$H^1_p(M)\simeq \{\alpha \in L^p(M; T^*M) \mid
d\alpha=0\ {\rm  and}\  d^*\alpha=0\}.$$
\end{proposition}
In general one cannot compare $L^p$ cohomology for different values of $p$.
However let us now assume that the Ricci curvature of $M$ is bounded from below:
\begin{equation}
\label{ricci}
\Ricci\ge -(n-1)\kappa^2 g. 
\end{equation}
Fro{}m the Bochner formula
$$\Delta \alpha \equiv \nabla^*\nabla \alpha+\Ricci (\alpha, \cdot) $$ we see that if $\alpha$ is a harmonic $1$-form in
$L^2$, then it satisfies the subelliptic estimate
\begin{equation} 
\Delta |\alpha|\le (n-1)\kappa^2 |\alpha| .
\label{subelliptic} 
\end{equation} 
Indeed if for $\varepsilon>0$ we define $f_\varepsilon= \sqrt{|\alpha|^2+\varepsilon}$, then it is classical to show that the Bochner formula and the Kato inequality imply
$$\Delta f_\varepsilon\le (n-1)\kappa^2 f_\varepsilon.$$ 
Passing to the limit $\varepsilon=0$, we get the desired subelliptic estimate (only in the distributional sense).
Hence with the Nash inequality we can deduce that $\alpha$ is in fact bounded and
\begin{equation}
\|\alpha\|_{L^\infty} \le C(n,\kappa,\mu) \|\alpha\|_{L^2}; 
\label{bounded}\end{equation} 
this can be done using a Nash-Moser iteration scheme \cite{Berard}, but 
we can also use our upper bound on the heat operator with the inequality (\ref{subelliptic}), to
assert that for every $x\in M$ $$t\mapsto \left(e^{-t(\Delta-(n-1)\kappa^2 )} 
|\alpha|\right)(x) $$ is non-decreasing. With the mapping properties of the heat operator (\ref{heatmapping}) we obtain
 $$|\alpha|(x)\le C t^{-n/2}e^{t(n-1)\kappa^2} 
\|\alpha\|_{L^2}; $$ 
with $1/t=(n-1)\kappa^2$, we obtain the desired bound. 
Hence there is a well defined map $$\cH^1(M)\rightarrow H^1_p(M).$$ Proposition~\ref{HdR} immediately implies

\begin{lemma}\label{inj}
Assume that (\ref{nash}, \ref{volumgro}, \ref{ricci}) hold, and that for some $p>2$, the Riesz transform is bounded on $L^p$ and on $L^{p/(p-1)}$. Then the natural map $$\cH^1(M)\rightarrow H^1_{p}(M)$$ is
injective.\end{lemma}

A corollary of this lemma is:

\begin{corollary}\label{nb}Assume  that $n>2$, and that $(M,g)$ satisfies the assumptions
(\ref{nash}, \ref{volumgro}, \ref{ricci}) and that it has more than two ends, then for
every $p\ge n$ the Riesz transform is not bounded on $L^p$.
\end{corollary}

 \proof In $n>2$, then $(M,g)$ satisfies the Sobolev inequality (\ref{sobo}),
  and following \cite{CSZ}, we know that if $M$ has more than two ends there
  exists a non-constant bounded harmonic function $h$ with finite Dirichlet
   energy\footnote{In fact if $M\setminus K=U_+\cup U_-$ with $K$ compact and
    $U_\pm$ unbounded then $\lim_{x\to\infty,x\in U_\pm} h(x)=\pm 1$.}, hence $dh$
     is a harmonic $L^2$ $1-$form. Take $\chi_k$ as in \eqref{cutoffk}.   We have $d\chi_kh=\chi_k dh+hd\chi_k$, but if $p\ge n$, then if $V(r):=\vol B(o,r)$ we
     have
     $$\| hd\chi_k\|^p_{L^p}\le C\|h\|^p_{L^\infty}  \int_k^{k^2} ((\log k)r)^{-p}dV(r)$$
     and integrating by parts we have :
     $$\int_k^{k^2} r^{-p}dV(r)= \frac{ \vol B(o,k^2)}{k^{2p}}-
     \frac{ \vol B(o,k)}{k^{p}}+p\int_k^{k^2}  \frac{ \vol B(o,r)}{r^{p+1}}dr.$$ This
     quantity is bounded for $p>n$ and grows as $C\log k$ if $p=n$. Hence we obtain that
$dh=L^p-\lim_{k\to\infty} d(\chi_kh)$.
Hence $dh$ is zero in reduced $L^p$ cohomology.  Since $dh$ is non-zero  in
$\cH^1(M)$, and since under (\ref{nash},\ref{volumgro})  the Riesz transform is bounded on $L^p$ for $1<p<2$, Lemma (\ref{inj}) says that it  can not be bounded on $L^p$ if $p\ge n$.\endproof
This generalizes the unboundedness part of Theorem~\ref{main} to the much larger class of manifolds satisfying \eqref{nash}, \eqref{volumgro}, \eqref{ricci}.

\section{Concluding remarks and open problems}\label{conc}In this final section we discuss some questions suggested
 by the work above, and pose some open problems. 
 
 For any complete Riemannian manifold $M$ of infinite measure there are numbers $\pmin \leq 2 \leq \pmax$ such
that the Riesz transform is bounded on $L^p$ for all $p$ between $\pmin$ and $\pmax$ (it may or may not be bounded at $p = \pmin$ or $p = \pmax$). We may call these values the lower and upper thresholds for $M$. 

There are a number of classes of manifolds on which the Riesz 
 transform is known to be bounded on $L^p$ for all $p$ (in other words, $\pmin = 1$ and $\pmax = \infty$); for example, manifolds with nonnegative Ricci curvature 
 \cite{Bakry},  Cartan-Hadamard manifolds with a spectral gap \cite{Lo}, noncompact symmetric spaces \cite{Anker}  and Lie 
 groups of polynomial growth \cite{Alex} (see \cite{ACDH} for more examples and references). 
On the other hand, Coulhon and Ledoux  showed in \cite{CL} that for any $p_0>2$ there is a manifold $M$ with bounded geometry such that $\pmax \leq p_0$.  Another example, with polynomial volume growth, was given in \cite{CDfull}.
     Since, as we mentioned earlier, it was shown by Coulhon and Duong that $\pmin = 1$ for a large class of
  complete manifolds, and that $\pmax = 2$ for certain simple surfaces, one could wonder whether $2$ is the  upper threshold for a large class of manifolds.  But H.-Q. Li \cite{Li} proved that for $n$-dimensional cones with compact basis, 
$$\pmax =  \begin{cases} n \Big( \frac{n}{2} - \sqrt{\big(\frac{n-2}{2} 
    \big)^2 + \lambda_1} \Big)^{-1}, \quad \lambda_1 < n-1 \\+\infty, \qquad \lambda_1 \geq n-1, \end{cases}$$where $\lambda_1$
     is the smallest nonzero eigenvalue of the Laplacian on the basis. Note that  $\pmax>n$ here.

     \begin{openproblem} Is a result similar to H.-Q. Li's valid  for smooth manifolds with one conic or asymptotically conic end? what happens for several conic ends? \end{openproblem}
     
 Manifolds with more than one Euclidean end satisfy the doubling condition but not the scaled $L^2$ Poincar\'e inequality. Our result sheds some light on the implications for the Riesz transform of these conditions: it follows from \cite{Li} and \cite{CouLi}  that doubling together with Poincar\'e (equivalently, upper and lower Gaussian estimates of the heat kernel) are not sufficient
     for the Riesz transform to be bounded for {\it all} $p>2$.   Theorem \ref{main} shows that 
     these conditions are not necessary for the Riesz transform to be bounded for {\it some} $p>2$.

     The class of manifolds with Euclidean ends is of course extremely special. One can  attempt to enlarge the class of known examples synthetically, i.e. by creating further examples from known examples by performing various operations.  Our results may be seen as obstructions to the stability of the $L^p$ boundedness of the Riesz transform under gluing for $p$ above the dimension.

     \begin{openproblem} Under which conditions is boundedness of the Riesz transform on $L^p$ stable 
     under the following operations on manifolds:
     \begin{itemize}
     \item gluing, 
     \item compact metric perturbations,
     \item taking products, $((M_1, g_1), (M_2, g_2)) \to (M_1 \times M_2, g_1 \oplus g_2)$,
     \item warped products.\end{itemize}\end{openproblem}

We only mention, without proof, one result  along these lines.
 Namely, if the Riesz transform is bounded on $L^p$ on a complete Riemannian manifold $M$ 
of infinite measure, and with Ricci curvature bounded from below,
then it is bounded on $L^p$ on $M \times N$ for any compact $N$.
   
We also mention a conjecture on manifolds obtained by gluing several copies of a 
simply connected nilpotent Lie group (endowed with a left invariant metric).
    According to \cite{Alex} we know that the Riesz transform is bounded 
     for every $p$ on a simply connected nilpotent Lie group.
     Let $(N,g_0)$ be a simply connected nilpotent Lie group of dimension $n>2$
      (endowed with a left invariant metric).  According to \cite{Alex} 
      we know that the Riesz transform on $(N,g_0)$ is bounded 
     for every $p$.  Let
     $\nu$ be the homogeneous dimension of $N$;  for instance we can set
     $$\nu=\lim_{R\to \infty} \frac{\log \vol B(o,R)}{\log R},$$
     $o\in N$ being a fixed point. Let $(M,g)$ be a manifold obtained by gluing $k>1$ copies of
     $(N,g_0)$. The manifold $M$ is diffeomorphic to the sphere with $k$ points removed.
     According to \cite{CD} we know that on $(M,g)$ the Riesz transform is bounded on $L^p$ for $p\in ]1,2]$. 
     The argument of Section~\ref{coh} can be applied (changing the $n$ appearing in analytic inequalities by the
     homogeneous dimension $\nu$). Then  (\ref{nb}) implies that the Riesz transform is not bounded 
     on $L^p$ for $p\ge \nu$. We can moreover compute the $L^p$ cohomology of $M$.
     The arguments of Proposition 3.3 of \cite{Canew}
     are given for the $L^2$ cohomology of such $M$ but they can be easily modified for $L^p$ cohomology.
We find that for $p\in ]1,\nu[$, $$H^1_p(M)\simeq H^1_c(M)=\R^{k-1},$$ whereas for 
     $p\ge \nu$ we find $$H^1_p(M)=\{0\}.$$
    That is, the conclusion of Lemma~\ref{inj} is satisfied if and only if $p<\nu$. This gives another proof
     of the fact that the Riesz transform is not bounded on $L^p$ for $p\ge \nu$. 

\begin{openproblem} Show that the Riesz transform on $(M,g)$ is bounded on $L^p$ for $p\in]1,\nu[$.
\end{openproblem}

Finally, it would be interesting to get results for differential forms. One can consider either $\nabla \Delta^{-1/2}$ or  $(d + \delta) \Delta^{-1/2}$, where $\Delta = (d + \delta)^2$ is the Laplacian on forms and $\Delta^{-1/2}$ really means $f(\Delta)$ where $f(0) = 0$ and $f(x) = x^{-1/2}$ for $x > 0$ (this projects off the $L^2$ null space of $\Delta$, i.e. the $L^2$-cohomology, which is trivial in the case of $0$-forms when $M$ has infinite measure).

      \begin{openproblem}Determine the upper and lower thresholds for $(d + \delta) \Delta^{-1/2}$ or $\nabla \Delta^{-1/2}$ acting on $k$-forms on manifolds with Euclidean ends.\end{openproblem}
      
      \begin{openproblem} Extend the results of section~\ref{coh} to differential forms of all degrees.\end{openproblem}

\end{document}